\documentclass[journal]{IEEEtran}

\usepackage{tikz}
\graphicspath{{images/}}
\usepackage{amsmath}
\usepackage{latexsym, amssymb}
\usepackage{graphicx}
\usepackage{amsmath, amsbsy}
\usepackage{amsopn, amstext}
\usepackage{ifpdf,hyperref}
\usepackage{cancel, color}
\usepackage{epstopdf}

\def\Char{Char}
\def\GL{GL}
\def\gl{gl}

\DeclareMathOperator{\Col}{Col}
\DeclareMathOperator{\Row}{Row}

\DeclareMathOperator{\lcm}{lcm}

\def\cal{\mathcal}

\def\ra{\rightarrow}

\def\a{\alpha}
\def\b{\beta}
\def\d{\delta}

\def\0{{\bf 0}}

\def\i{{\bf i}}
\def\j{{\bf j}}
\newcommand{\R}{{\mathbb R}}

\newcommand{\C}{{\mathbb C}}
\newcommand{\F}{{\mathbb F}}
\def\H{\mbox{\boldmath$H$}}

\def\dsum{\mathop{\sum}\limits}

\newtheorem{thm}{Theorem}[section]
\newtheorem{dfn}[thm]{Definition}
\newtheorem{prp}[thm]{Proposition}
\newtheorem{exa}[thm]{Example}

\newtheorem{rem}[thm]{Remark}

\begin{document}

\title{On Perfect Hypercomplex Algebra}

\author{Daizhan Cheng,  Zhengping Ji
	\thanks{This work is supported partly by the National Natural Science Foundation of China (NSFC) under Grants 62073315, 61074114, and 61273013.}
	\thanks{Key Laboratory of Systems and Control, Academy of Mathematics and Systems Sciences, Chinese Academy of Sciences,
		Beijing 100190, P. R. China (e-mail: dcheng@iss.ac.cn, jizhengping@amss.ac.cn).}
}


\maketitle

\begin{abstract}
 The set of associative and commutative hypercomplex numbers, called the perfect hypercomplex algebra (PHA) is investigated. Necessary and sufficient conditions for an algebra to be a PHA  via semi-tensor product(STP) of matrices are reviewed. The zero set is defined for non-invertible hypercomplex numbers in a given PHA, and a characteristic function is proposed for calculating zero set. Then PHA of different dimensions are considered. First, $2$-dimensional PHAs are considered as examples to calculate their zero sets etc. Second, all the $3$-dimensional PHAs are obtained and the corresponding zero sets are investigated. Third, $4$-dimensional or even higher dimensional PHAs are also considered. Finally, matrices over pre-assigned PHA, called perfect hypercomplex matrices (PHMs) are considered. Their properties are also investigated.
\end{abstract}


\IEEEpeerreviewmaketitle

\section{Introduction}

Hypercomplex number (HN) are generalization of complex numbers ($\C$). A class of HN with pre-assigned addition and product form a special vector space over $\R$, called a hypercomplex algebra (HA).  HAs have various applications including signal and image processing \cite{pei04},  dealing with differential operators \cite{ala02,abr15}, designing neural networks \cite{cas19}, etc.

It was proved by Weierstrass  that the only finite field extension of real numbers ($\R$) is complex numbers ($\C$) \cite{li99}.  HA can be considered as an extension of real numbers ($\R$) to finite dimensional algebras. We call such extension finite algebra extension of real numbers.

In this paper we consider only a particular kind of finite algebra extensions of $\R$, which are commutative and associative. Hence, throughout this paper the following is assumed:

\vskip 2mm

{\bf Assumption 1}: $\H$ is the set of PHAs, that is, the set of finite dimensional algebras over $\R$ which are commutative and associative.

\vskip 2mm

In addition to complex numbers, hyperbolic numbers, dual numbers, and Tessarine quaternion are also PHA.

STP of matrices is a generalization of conventional matrix product. It is a powerful tool to deal with multi-linear mappings. In \cite{che21}, STP was used to investigate finite algebra extensions of $\R$. (in fact, the extensions of any $\F$ with $\Char(F)=0$ have been discussed there.) In \cite{fu21}, STP was used to investigate general Boolean-type algebras. A key issue in these approaches is to define a matrix, called product matrix of a certain algebra. Using STP, associativity, commutativity, and some other properties of a finite algebra extension can be verified via its product matrix.

In this paper, this STP approach is used to investigate PHA and PHM. For PHA, first the formulas for verifying whether an HA is associative and commutative are reviewed. Then the zero set is defined as the set of non-invertible numbers. A characteristic function is proposed to calculate (or characterize) the zero set. Then the PHA of dimensions 2, 3, or 4 are constructed separately. Even higher dimensional cases are also discussed. Their zero-sets, which are of measure zero, are calculated. Analytic functions and some other properties of PHAs are then discussed.

As for PHMs, their invertibility and some further properties are investigated. The Lie-group and Lie algebra over a pre-assigned PHA, say, ${\cal A}$, are also defined as general linear group, $\GL(n,{\cal A})$ and general linear algebra, $\gl(n,{\cal A})$ are investigated.  Certain properties are revealed.

The rest of this paper is organized as follows:

Before ending this section, we give a list of notations:

\begin{enumerate}

\item $\H$: Set of PHAs, as finite algebra extensions over $\R$.

\item ${\cal A}_{m\times n}$: set of $m\times n$ dimensional matrices, with all entries in algebra ${\cal A}$.

\item $\ltimes$: STP of matrices.

\item $\Col(A)$ : the set of columns (rows) of ~$A$; $\Col_i(A)$ ($\Row_i(A)$): the $i$-th column (row) of ~$A$.


%
\item
$\d_k^i$: The $i$-th column of identity matrix $I_k$.


%
%
%
%

\end{enumerate}

\vskip 2mm

\section{Preliminaries} 

\subsection{Semi-tensor Product of Matrices}

Since STP is a fundamental tool in this approach, this section will give a brief survey for STP. We refer to \cite{che12} for more details.
 STP is a generalization of conventional matrix product, defined as follows:

\begin{dfn}\label{d2.1.1} Let $A\in \R_{m\times n}$ and $B\in \R_{p\times q}$, $t=\lcm(n,p)$ be the least common multiple of $n$ and $p$. Then the STP of $A$ and $B$, denoted by $A\ltimes B$, is defined by
\begin{align}\label{2.1.1}
A\ltimes B:=\left(A\otimes I_{t/n}\right)\left(B\otimes I_{t/P}\right),
\end{align}
where $\otimes$ is Kronecker product.
\end{dfn}

It is easy to see that STP is a generalization of conventional matrix product. That is, when $n=p$, the STP is degenerated to the conventional matrix product, i.e.,, $A\ltimes B=AB$. Because of this, in most cases the symbol $\ltimes$ is omitted.

One of the most important advantages of STP is that STP keeps most important properties of conventional matrix product available, including association, distribution, etc. In the following we introduce some additional properties of STP, which will be used in the sequel.

Define a swap matrix $W_{[m,n]}\in {\cal M}_{mn\times mn}$ as follows:
\begin{align}\label{2.1.2}
W_{[m,n]}:=\left[I_n\otimes \d_m^1, I_n\otimes \d_m^2,\cdots, I_n\otimes \d_m^m,\right]\in {\cal L}_{mn\times mn}.
\end{align}

\begin{prp}\label{p2.1.2} Let $x\in \R^m$ and $y\in \R^n$ be two column vectors. Then
\begin{align}\label{2.1.3}
W_{[m,n]}x\ltimes y=y\ltimes x.
\end{align}
\end{prp}

The following proposition ``swaps" a vector with a matrix:

\begin{prp}\label{p2.1.3} Let $x\in \R^t$  be a column vector, and $A$ be an arbitrary matrix. Then
\begin{align}\label{2.1.4}
x\ltimes A=(I_t\otimes A)\ltimes x.
\end{align}
\end{prp}

Throughout this paper the default matrix product is assumed to be STP, and the symbol $\ltimes$ is omitted if there is no possible confusion.

\subsection{Matrix Expression of an Algebra}

We are only interested in algebras over $\R$.

\begin{dfn}\label{d2.2.1}\cite{hun74}
\begin{itemize}
\item[(i)] An algebra over $\R$ is a pair, denoted by ${\cal A}=(V,*)$,  where  $V$ be a real vector space, and $*:V\times V\ra V$, satisfies
\begin{align}\label{2.2.1}
\begin{array}{l}
(ax+by)*z=ax*z+by*z,\\
x*(ay+bz)=ax*y+bx*z,\quad x,y,z\in V,\;a,b\in \R.\\
\end{array}
\end{align}

\item[(ii)] An algebra ${\cal A}=(V,*)$ is said to be commutative, if
\begin{align}\label{2.2.2}
x*y=y*x,\quad x,y\in V.
\end{align}
\item[(iii)] An algebra ${\cal A}=(V,*)$ is said to be associative, if
\begin{align}\label{2.2.3}
(x*y)*z=x*(y*z),\quad x,y,z\in V.
\end{align}
\end{itemize}
\end{dfn}

\begin{dfn}\label{d2.2.2} Let ${\cal A}=(V,*)$ be a $k$-dimensional vector space with $e=\{\i_1,\i_2,\cdots,\i_k\}$ as a set of basis. Denote
\begin{align}\label{2.2.4}
\i_i*\i_j=\dsum_{s=1}^kc^s_{i,j}\i_s,\quad i,j=1,2,\cdots,k.
\end{align}
Then the product matrix of ${\cal A}$ is defined as
\begin{align}\label{2.2.5}
P_{{\cal A}}:=
\begin{bmatrix}
c^1_{1,1}&c^1_{1,2}&\cdots&c^1_{1,k}&\cdots&c^1_{k,k}\\
c^2_{1,1}&c^2_{1,2}&\cdots&c^2_{1,k}&\cdots&c^2_{k,k}\\
~&~&\ddots&~&\ddots&~\\
c^k_{1,1}&c^k_{1,2}&\cdots&c^k_{1,k}&\cdots&c^k_{k,k}\\
\end{bmatrix}.
\end{align}
\end{dfn}

Assume $x=\dsum_{j=1}^kx_i\i_j$, is expressed in a column vector form as $x=(x_1,x_2,\cdots,x_k)^T$. Similarly, $y=(y_1,y_2,\cdots,y_k)^T$. Then

\begin{thm} \label{t2.2.3}\label{che21} In vector form two product of two hypercomplex numbers $x,y\in {\cal A}$ is computable via following formula. 

\begin{align}\label{2.2.6}
x*y=P_{{\cal A}}xy.
\end{align}
\end{thm}

Using formula (\ref{2.2.6}) and the properties of STP yields the following results, which are fundamental for our further investigation.

\begin{thm} \label{t2.2.4}\label{che21} 
\begin{itemize}
\item[(i)] ${\cal A}$ is commutative, if and only if,
\begin{align}\label{2.2.7}
P_{{\cal A}}\left[I_k-W_{[k,k]}\right]=0.
\end{align}
\item[(ii)]  ${\cal A}$ is associative, if and only if,
\begin{align}\label{2.2.8}
P^2_{{\cal A}}=P_{{\cal A}}\left(I_k \otimes P_{{\cal A}}\right).
\end{align}
\end{itemize}
\end{thm}

\section{Hypercomplex Numbers}

\subsection{Perfect Hypercomplex Algebra on $\R$}

\begin{dfn}\label{d3.1.1} \cite{she89} 
A number $p$ is called a hypercomplex number, if it can be expressed in the form
\begin{align}\label{3.1.1}
p=p_0+p_1\i_1+\cdots+p_n\i_n,
\end{align}
where $p_i\in \R$, $i=0,1,\cdots,n$, $\i_i$, $i=1,2,\cdots,n$ are called hyperimaginary units.
\end{dfn}

\begin{rem}\label{r3.1.2} A hypercomplex number may belong to different algebras, depending on their product structure matrices. A hypercomplex algebra, denoted by ${\cal A}$, is an algebra over $\R$ with basis $e=\{\i_0:=1,\i_1,\cdots,\i_n\}$. 
\end{rem}

\begin{prp}\label{p3.1.3} Assume 
$${\cal A}=\{p_0+p_1\i_1+\cdots+p_n\i_n\;|\;p_0,p_1,\cdots,p_n\in \R\}.
$$
Then its product matrix 
$$
P_{{\cal A}}:=[M_0,M_1,\cdots,M_n],
$$
where $M_i\in \R_{(n+1)\times (n+1)}$, $i=0,1,\cdots,n$, satisfies
 the following condition.
\begin{itemize}
\item[(i)] 
\begin{align}\label{3.1.2}
M_0=I_{n+1}
\end{align}
is an identity matrix.
\item[(ii)]
\begin{align}\label{3.1.3}
\Col_1(M_j)=\d^{j+1}_{n+1},\quad j=1,2,\cdots,n.
\end{align}
\end{itemize}
\end{prp}

\begin{dfn}\label{d3.1.4} A hypercomplex algebra ${\cal A}$ is called a PHA, denoted by ${\cal A}\in \H$, if it is commutative and associative.
\end{dfn}

\begin{exa}\label{e3.1.5} Consider $\C$. It is easy to calculate that its product structure matrix is
\begin{align}\label{3.1.4}
P_{\C}=\begin{bmatrix}
1&0&0&-1\\
0&1&1&0
\end{bmatrix}.
\end{align}
A straightforward computation verifies (\ref{2.2.7}) and (\ref{2.2.8}), hence it is a PHA.
\end{exa}

\subsection{Invertibility of elements in PHA}

Now for a PHA, say, ${\cal A}=(V,*)$, if every $0\neq x\in V$ has its inverse $x^{-1}$ such that $x*x^{-1}=x^{-1}x=1$, then ${\cal A}$ is a field. Unfortunately, according to Weierstrass, if ${\cal A}\neq \C$, it is not a field. Hence, when an element $x\in V$ (we also say $x\in {\cal A}$, which means $x\in V$.) is invertible. 

To answer this question, we need some new concepts, which are firstly discussed in  \cite{che21}.

\begin{dfn}\label{d3.2.1} 
\begin{itemize}
\item[(i)] Let $A_1,A_2,\cdots,A_r$ be a set of square real matrices.  $A_1,A_2,\cdots,A_r$ are said to be jointly nonsingular, if their non-trivial linear combination is non-singular. That is, assume
$$
det\left(\dsum_{i=1}^rc_iA_i\right)=0,
$$ 
then $c_1=c_2=\cdots=c_r=0$.
\item[(ii)] Let $A\in \R_{k\times k^2}$. $A$ is said to be jointly nonsingular, if 
$A=[A_1,A_2,\cdots,A_k$, where $A_s\in \R_{k\times k^2}$ and $\{A_i\|i=1,2,\cdots,k\}$ are jointly nonsingular.
\end{itemize}
\end{dfn}

Then we have the following result:

\begin{prp}\label{p3.2.2}\cite{che21} ~$A\in \R_{k\times k^2}$ is jointly non-singular, if and only if, 
one of the following two equivalent conditions is satisfied:
\begin{itemize}
\item[(i)] The matrix $P_{{\cal A}}x\in \R_{k\times k}$ is non-singular for all $x\neq 0$.
\item[(ii)]
The following homogeneous polynomial
\begin{align}\label{3.2.1}
\begin{array}{l}
\xi(x_1,\cdots,x_k)=\det(Ax)\\
~~=
\dsum_{i_1=1}^k\cdots \dsum_{i_k=1}^k\mu_{i_1,\cdots,i_k}x_{i_1}\cdots x_{i_k}\neq 0,\quad \forall x\neq 0.
\end{array}
\end{align}
\end{itemize}
\end{prp}

The $\xi(x_0,x_1,\cdots,x_n)$ is called the characteristic function of ${\cal A}$.

\begin{exa}\label{e3.2.3}
Consider  ~$\C=\R(i)$. ~~Calculating right hand side of (\ref{3.2.1}) for $P_{\C}$, we have
$$
xi(x_1,x_2)=x_1^2+x_2^2.
$$
Hence, $\xi(x_1,x_2)=0$, if and only if, $x_1=x_2=0$. It follows that $P_{\C}$ is jointly non-singular.
\end{exa}

Summarizing above arguments, we have the following result.

\begin{prp}\label{p3.2.4} Let ${\cal A}$ be a finite dimensional algebra over $\R$. Then ${\cal A}$ is a field, if and only if,
\begin{itemize}
\item[(i)] ${\cal A}$ is commutative, that is, (\ref{2.2.7}) holds;
\item[(ii)]  ${\cal A}$ is associative, that is, (\ref{2.2.8}) holds;
\item[(iii)] Each $0\neq x\in {\cal A}$ is invertible, that is, $P_{{\cal A}}$ is jointly invertible.
\end{itemize}
\end{prp}

Unfortunately, it is well known that the only finite dimensional algebra over $\R$, which is a field, is $\C$. In this paper we are particularly interested in commutative and associative algebra ${\cal A}\in \H$. Then unless ${cal A}=\C$, there must be some elements $x\neq 0$, which are not invertible.

\begin{dfn}\label{d3.2.5} Let ${\cal A}\in \H$. Its zero set is defined by
\begin{align}\label{3.2.2}
{\cal Z}_{{\cal A}}:=\left\{z\in {\cal A}\;|\; \det\left(P_{{\cal A}}z\right)=0\right\}
\end{align}
\end{dfn}

It is clear that
\begin{itemize}
\item[(i)]
if ${\cal A}=\C$, then ${\cal Z}_{{\cal A}}=\{0\}$;
\item[(ii)] if  ${\cal A}\neq \C$, then $\{0\}\in {\cal Z}_{{\cal A}}$, and
${\cal Z}_{{\cal A}}\backslash\{0\}\neq \emptyset$.
\end{itemize}

\subsection{Coordinate Transformation of Hypercomplex Algebra}

Since the product matrix $P_{{\cal A}}$ of a hypercomplex algebra depends on the basis, it is necessary to consider different forms of $P_{{\cal A}}$ under a change of basis, which is commonly called a coordinate transformation. Let $\i_1,\i_2,\cdots,\i_n$ be the hyperimaginary units of ${\cal A}$ and
\begin{align}\label{3.3.1}
(1,\j_1,\j_2,\cdots,\j_n)=(1,\i_1,\i_2,\cdots,\i_n)T,
\end{align}
\begin{align}\label{3.3.2}
T=\begin{bmatrix}
1&E\\
0&T_0
\end{bmatrix}
\in \R_{(n+1)\times (n+1)},
\end{align}
where $T_0$ is non-singular.

Then $x\in {\cal A}$ can be expressed as
$$
x=x_0+x_1\i_1+\cdots+x_n\i_n=x_0+\bar{x}_1\j_1+\cdots+\bar{x}_n\j_n,
$$
which is called a coordinate change on ${\cal A}$.

Denote $x=(x_0,x_1,\cdots,x_n)^T$, $\bar{x}=(\bar{x}_0,\bar{x}_1,\cdots,\bar{x}_n)^T$. Then
\begin{align}\label{3.3.3}
\bar{x}=T^{-1}x.
\end{align}

\begin{dfn}\label{d3.3.1} Let ${\cal A}$ and $\overline{{\cal A}}$ be two $n+1$ dimensional hypercomplex algebras. ${\cal A}$ and $\overline{{\cal A}}$  are called isomorphic, if there exists a bijective mapping $\Psi:{\cal A} \ra \overline{{\cal A}}$, satisfying
\begin{itemize}
\item[(i)]
\begin{align}\label{3.3.4}
\Psi(1)=1
\end{align}
\item[(ii)]
\begin{align}\label{3.3.5}
\Psi(ax+by)=a\Psi(x)+b\Psi(y),\quad x,y\in {\cal A}, a,b\in \R;
\end{align}
\item[(iii)]
\begin{align}\label{3.3.6}
\Psi(x*y)=\Psi(x)*\Psi(y),\quad x,y\in {\cal A}.
\end{align}
\end{itemize}
$\Psi$ is called an isomorphism.
\end{dfn}

A straightforward verification shows the following result immediately.

\begin{prp}\label{p3.3.2} 
 Assume
$$
\begin{array}{l}
{\cal A}=\left\{x_0+\dsum_{i=1}^nx_i\i_i\;|\;x_0,x_1£¬\cdots, x_n\in \R\right\}\\
\overline{{\cal A}}=\left\{\bar{x}_0+\dsum_{i=1}^n\bar{x}_i\j_i\;|\;\bar{x}_0,\bar{x}_1£¬\cdots, \bar{x}_n\in \R\right\}.\\
\end{array}
$$
${\cal A}$ and $\overline{{\cal A}}$ are isomorphic, if and only if, there is a non-singular matrix $T\in \R_{(n+1)\times (n+1)}$ of the form (\ref{3.3.2}), such that (\ref{3.3.3}) holds true.
\end{prp}

\begin{prp}\label{p3.3.3} Assume ${\cal A},~\overline{{\cal A}}\in \H_{n+1}$, with their PSMs as
$P_{{\cal A}}$ and $P_{\overline{{\cal A}}}$ respectively. ${\cal A}$ and $\overline{{\cal A}}$ are isomorphic, if and only if, there exists a non-singular matrix
$T$ as in (\ref{3.3.2})
such that
\begin{align}\label{3.3.7}
P_{\overline{{\cal A}}}=T^{-1}P_{{\cal A}}\left(T\otimes T\right).
\end{align}
\end{prp}

\begin{IEEEproof}
(Necessary) Let $T$ be constructed as is (\ref{3.3.2}), such that
$$
\bar{x}=T^{-1}x.
$$
Then we have
\begin{align}\label{3.3.8}
P_{{\cal A}}xy=TP_{\overline{\cal A}}\bar x\bar y,\quad x,y\in {\cal A}.
\end{align}
The right hand side (RHS) of (\ref{3.3.8}) becomes
$$
\begin{array}{ccl}
RHS_{(\ref{3.3.8})}&=&=TP_{\overline{\cal A}}T^{-1}xT^{-1}y\\
~&=&TP_{\overline{\cal A}}T^{-1}\left(I_{n+1}\otimes T^{-1}\right)xy.
\end{array}
$$
Since $x,y$ are arbitrary, we have
$$
P_{{\cal A}}=TP_{\overline{\cal A}}T^{-1}\left(I_{n+1}\otimes T^{-1}\right).
$$
Hence,
$$
\begin{array}{ccl}
P_{\overline{\cal A}}&=& T^{-1} P_{{\cal A}}\left(I_{n+1}\otimes T\right)T\\
~&=& T^{-1} P_{{\cal A}}\left(T\otimes T\right)\\
\end{array}
$$

(Sufficiency) If (\ref{3.3.8}) holds true, it is easy to verify that
$$
\bar{x}=T^{-1}x
$$
is an isomorphism.

\end{IEEEproof}

\section{Lower Dimensional PHAs}

This section considers some examples of various dimensional algebras. 

\subsection{Structure of ${\cal A}\in \H_2$}

Consider ${\cal A}\in \H_2$. According to Proposition\ref{p3.1.3}, its PSM is
\begin{align}\label{4.1.1}
P_{{\cal A}}=\begin{bmatrix}
1&0&0&\a\\
0&1&1&\b
\end{bmatrix}
\end{align}

Consider a coordinate change
$$
T=\begin{bmatrix}
1&s\\
0&t\end{bmatrix},\quad t\neq 0.
$$
Using formula (\ref{3.3.7}), we have
\begin{align}\label{4.1.2}
\begin{array}{ccl}
P_{\overline{\cal A}}&=&
T^{-1}P_{{\cal A}}\left(T\otimes T\right)\\
~&=&\begin{bmatrix}
1&0&0&\a t^2-s(s+t\b)\\
0&1&1&s+t\b\\
\end{bmatrix}
\end{array}
\end{align}
Now we may choose
$$
s=-t\b,
$$
Then we have
\begin{align}\label{4.1.3}
P_{\overline{\cal A}}=
\begin{bmatrix}
1&0&0&\a t^2\\
0&1&1&0\\
\end{bmatrix}
\end{align}

We consider them case by case.

\begin{itemize}

\item If $\a=0$, we have
\begin{align}\label{4.1.4}
P_{\overline{\cal A}}=
\begin{bmatrix}
1&0&0&0\\
0&1&1&0\\
\end{bmatrix}
\end{align}

\item If $\a>0$, choosing $t=\frac{1}{\sqrt{|\a|}}$, then we have that
\begin{align}\label{4.1.5}
P_{\overline{\cal A}}=
\begin{bmatrix}
1&0&0&1\\
0&1&1&0\\
\end{bmatrix}.
\end{align}

\item If $\a<0$, choosing $t=\frac{1}{\sqrt{|\a|}}$ yields
\begin{align}\label{4.1.6}
P_{\overline{\cal A}}=
\begin{bmatrix}
1&0&0&-1\\
0&1&1&0\\
\end{bmatrix}.
\end{align}
\end{itemize}

We conclude that up to isomorphism there are three ${\cal A}\in \H_2$, they are
\begin{itemize}
\item set of dual numbers (${\cal A}_D$), which corresponds to (\ref{4.1.4});
\item set of hyperbolic numbers (${\cal A}_H$), which corresponds to (\ref{4.1.5});
\item set of complex numbers ($\C$), which corresponds to (\ref{4.1.6}).
\end{itemize}

Next, using (\ref{3.2.1}), we can calculate their characteristic functions.
\begin{itemize}
\item
\begin{align}\label{4.1.7}
\xi_{{\cal A}_D}=x_0^2.
\end{align}
Then
\begin{align}\label{4.1.8}
{\cal Z}_{{\cal A}_D}=\{x_0+x_1\i\in {\cal A}_D\;|\; x_0=0\}.
\end{align}

\item
\begin{align}\label{4.1.9}
\xi_{{\cal A}_H}=x_0^2-x_1^2.
\end{align}
Then
\begin{align}\label{4.1.10}
{\cal Z}_{{\cal A}_D}=\{x_0+x_1\i\in {\cal A}_H\;|\; x_0=\pm x_1\}.
\end{align}

\item
\begin{align}\label{4.1.11}
\xi_{\C}=x_0^2+x_1^2;.
\end{align}
Then
\begin{align}\label{4.1.12}
{\cal Z}_{\C}=\{0\}.
\end{align}

\end{itemize}

\begin{rem}\label{r4.1.1} 
\begin{itemize}
\item[(i)] It is obvious that all ${\cal A}_D$, ${\cal A}_H$, and $\C$ are associative and commutative, that is they are all PHAs.
\item[(ii)] They have minimum polynomials $x_0^2$, $x_0^2-x_1^2$, and $x_0^2+x_1^2$ respectively. According to Galois theory, only the minimum polynomial of $\i$ is irreducible, the extension is a field. Hence, only $\C$ is a fields.
\item[(iii)] It is easily seen that their zero sets are all zero measure set. This is always true for all PHA. Because their zero set are zero set of algebraic equations, which are called algebraic numbers. The set of algebraic numbers is always zero measure set.
\end{itemize}
\end{rem}

\subsection{Structure of Triternions}

\begin{dfn}\label{d4.2.1} An algebra  ${\cal A}$ of dimension $3$ is called a triternion if ${\cal A}\in \H_3$.
\end{dfn}

It is easy to see that if ${\cal A}\in \H_3$ is symmetric, then its PSM is
\begin{align}\label{4.2.1}
P_{{\cal A}}=
\begin{bmatrix}
1&0&0&0&a&d&0&d&p\\
0&1&0&1&b&e&0&e&q\\
0&0&1&0&c&f&1&f&r\\
\end{bmatrix}
\end{align}

Next, we consider when ${\cal A}$ is associative. According to Theorem \ref{t2.2.4},
the necessary and sufficient condition is
\begin{align}\label{4.2.2}
P^2_{{\cal A}}=P_{{\cal A}}\left(I_3\otimes P_{{\cal A}} \right).
\end{align}
Denote $I=I_3$,
$$
A=\begin{bmatrix}
0&a&d\\
1&b&e\\
0&c&f
\end{bmatrix},\quad
B=\begin{bmatrix}
0&d&p\\
0&e&q\\
1&f&r
\end{bmatrix}.
$$
A direct computation shows that
\begin{align}\label{4.2.3}
\begin{array}{l}
\mbox{LHS of (\ref{4.2.2})}=(I,A,B,A,aI+bA+cB,\\
dI+eA+fB,B,dI+eA+fB,pI+qA+rB),\\
\mbox{RHS of (\ref{4.2.2})}=(I,A,B,A,A^2,AB,B,BA,B^2).
\end{array}
\end{align}

Then we have the following result:

\begin{thm}\label{t4.2.2}
${\cal A}\in \H_3$, if and only if, $P_{{\cal A}}$ has the form of (\ref{4.2.1}) with parameters satisfying
\begin{align}\label{4.2.4}
\begin{array}{l}
a=ce+f^2-bf-cr,\\
d=cq-ef,\\
p=e^2+fq-bq-er.\\
\end{array}
\end{align}
\end{thm}

\begin{IEEEproof}
(Necessity) Comparing both sides of (\ref{5.2}), (\ref{5.3}) shows that a necessary condition for (\ref{4.2.2}) holds true is (refer to the 6th and 8th blocks of both sides)
\begin{align}\label{4.2.5}
AB=BA.
\end{align}
Then it is easy to verify that (\ref{4.2.4}) provides necessary and sufficient condition for (\ref{4.2.5}) to be true.

(Sufficiency) A careful computation shows as long as (\ref{4.2.4}) holds, the RHS of (\ref{4.2.2}) and the LHS of (\ref{4.2.2}), shown in (\ref{4.2.3}), are equal.

\end{IEEEproof}

\begin{rem}\label{r4.2.3} Theorem \ref{t4.2.2}) provides an easy way to construct ${\cal A}\in \H_3$.
In fact, the parameters $b,c,e,f,q,r$ can be assigned freely, then $a,d,p$ can be obtained by (\ref{4.2.4}).
It is easy to see that there are uncountably many algebras of dimension 3, which are commutative and associative.
\end{rem}

Next, we give a numerical example.

\begin{exa}\label{e4.2.4} Construct ${\cal A}\in \H_3$ by setting $b=c=f=q=r=0$ and $e=1$. Then we have
$d=a=0$ and $p=1$. The PSM of ${\cal A}$ is
\begin{align}\label{4.2.6}
P_{{\cal A}}=
\begin{bmatrix}
1&0&0&0&0&0&0&0&1\\
0&1&0&1&0&1&0&1&0\\
0&0&1&0&0&0&1&0&0\\
\end{bmatrix}.
\end{align}
In fact, when $x\in {\cal A}$ is expressed into standard form as
$$
x=x_0+x_1\i_1+x_2\i_2,\quad x_0,x_1,x_2\in \R,
$$
then we have
$$
\begin{array}{l}
\i_1^2=0;\quad \i_2^2=1,\\
\i_1*\i_2=\i_2*\i_1=\i_1.
\end{array}
$$

Then it is easy to calculate that
\begin{align}\label{4.2.7}
\xi_{{\cal A}}=(x_0-x_2)(x_0+x_2)^2.
\end{align}
Hence,
\begin{align}\label{4.2.8}
{\cal Z}_{{\cal A}}=\{(x_0,x_1,x_2)\in \R^3\;|\;x_0=\pm x_2\}.
\end{align}
\end{exa}

\subsection{Structure of Perfect Quaternions}

This section considers some some algebras in $\H_4$. It seems not easy to provide a general description for algebras in $\H_4$. The principle argument is similar to triternions. We give some simple examples.

\begin{exa}\label{e4.3.1}
Consider an ${\cal A}\in \H_4$. Assume
$$
{\cal A}=\left\{ p_0+p_1\i_1+p_2\i_2+p_3\i_3\;|\;p_0,p_1,p_2,p_3\in \R\right\},
$$
satisfying
$$
\begin{array}{ll}
\i_1^2,\i_2^2,\i_3^2\in \{-1,0,1\},&
\i_1*\i_2=\i_2*\i_1=\pm \i_3,\\
\i_2*\i_3=\i_3*\i_2=\pm \i_1,&
\i_3*\i_1=\i_1*\i_3=\pm \i_2.\\
\end{array}
$$
To save space, we denote
$$
P_{{\cal A}_i}=[I_4,Q_i].
$$
Using MATLAB for an exhausting searching,
we got eight algebras as follows:
\begin{itemize}
\item
{\tiny
$$
Q_1=\left[\begin{array}{cccccccccccc}
 0&-1& 0& 0& 0& 0&-1& 0& 0& 0& 0& 1\\
 1& 0& 0& 0& 0& 0& 0& 1& 0& 0& 1& 0\\
 0& 0& 0& 1& 1& 0& 0& 0& 0& 1& 0& 0\\
 0& 0&-1& 0& 0&-1& 0& 0& 1& 0& 0& 0
 \end{array}\right]
$$
}
\item
{\tiny
$$
Q_2=\left[\begin{array}{cccccccccccc}
 0&-1& 0& 0& 0& 0&-1& 0& 0& 0& 0& 1\\
 1& 0& 0& 0& 0& 0& 0&-1& 0& 0&-1& 0\\
 0& 0& 0&-1& 1& 0& 0& 0& 0&-1& 0& 0\\
 0& 0& 1& 0& 0& 1& 0& 0& 1& 0& 0& 0
 \end{array}\right]
$$
}
\item
{\tiny
$$
Q_3=\left[\begin{array}{cccccccccccc}
 0&-1& 0& 0& 0& 0&1& 0& 0& 0& 0&-1\\
 1& 0& 0& 0& 0& 0& 0&-1& 0& 0&-1& 0\\
 0& 0& 0& 1& 1& 0& 0& 0& 0& 1& 0& 0\\
 0& 0&-1& 0& 0&-1& 0& 0& 1& 0& 0& 0
 \end{array}\right]
$$
}
\item
{\tiny
$$
Q_4=\left[\begin{array}{cccccccccccc}
 0&-1& 0& 0& 0& 0&1& 0& 0& 0& 0&-1\\
 1& 0& 0& 0& 0& 0& 0&1& 0& 0&1& 0\\
 0& 0& 0&-1& 1& 0& 0& 0& 0&-1& 0& 0\\
 0& 0& 1& 0& 0& 1& 0& 0& 1& 0& 0& 0
 \end{array}\right]
$$
}

\item
{\tiny
$$
Q_5=\left[\begin{array}{cccccccccccc}
 0& 1& 0& 0& 0& 0&-1& 0& 0& 0& 0&-1\\
 1& 0& 0& 0& 0& 0& 0& 1& 0& 0&1& 0\\
 0& 0& 0&-1& 1& 0& 0& 0& 0&-1& 0& 0\\
 0& 0&-1& 0& 0&-1& 0& 0& 1& 0& 0& 0
 \end{array}\right]
$$
}

\item
{\tiny
$$
Q_6=\left[\begin{array}{cccccccccccc}
 0& 1& 0& 0& 0& 0&-1& 0& 0& 0& 0&-1\\
 1& 0& 0& 0& 0& 0& 0&-1& 0& 0&-1& 0\\
 0& 0& 0& 1& 1& 0& 0& 0& 0& 1& 0& 0\\
 0& 0& 1& 0& 0& 1& 0& 0& 1& 0& 0& 0
 \end{array}\right]
$$
}

\item
{\tiny
$$
Q_7=\left[\begin{array}{cccccccccccc}
 0& 1& 0& 0& 0& 0& 1& 0& 0& 0& 0& 1\\
 1& 0& 0& 0& 0& 0& 0&-1& 0& 0&-1& 0\\
 0& 0& 0&-1& 1& 0& 0& 0& 0&-1& 0& 0\\
 0& 0&-1& 0& 0&-1& 0& 0& 1& 0& 0& 0
 \end{array}\right]
$$
}

\item
{\tiny
$$
Q_8=\left[\begin{array}{cccccccccccc}
 0& 1& 0& 0& 0& 0& 1& 0& 0& 0& 0& 1\\
 1& 0& 0& 0& 0& 0& 0& 1& 0& 0& 1& 0\\
 0& 0& 0& 1& 1& 0& 0& 0& 0& 1& 0& 0\\
 0& 0& 1& 0& 0& 1& 0& 0& 1& 0& 0& 0
 \end{array}\right]
$$
}
\end{itemize}
\end{exa}

Next, choose some ${\cal A}\in \H_4$ for further study.

\begin{exa}\label{e4.3.2} Recall Example \ref{e4.3.1}).

\begin{itemize}

\item[(i)] Consider ${\cal A}_3$:

It is easy to calculate that
\begin{align}\label{4.3.1}
\begin{array}{ccl}
\xi_{{\cal A}_3}&=&\det(P_{{\cal A}_3}x)\\
~&=&(x_0^2-x_2^2)^2+(x_1^2-x_3^2)^2\\
~~&~&+2(x_0x_1+x_2x_3)^2+2(x_0x_3+x_1x_2)^2.
\end{array}
\end{align}
It follows that
\begin{align}\label{4.3.2}
\begin{array}{l}
{\cal Z}_{{\cal A}_3}=\left\{(x_0,x_1,x_2,x_3)^T\in \R^4\;|\;\right.\\
\left. (x_0=x_2)\cap(x_1=-x_3)~\mbox{or}~(x_0=-x_2)\cap(x_1=x_3) \right\}.
\end{array}
\end{align}

\item[(ii)] Consider ${\cal A}_8$:

It is easy to calculate that
\begin{align}\label{4.3.3}
\begin{array}{ccl}
\xi_{{\cal A}_8}&=&\det(P_{{\cal A}_8}x)\\
~&=&x_0^4+x_1^4+x_2^4+x_3^4\\
~~&~&-2(x^2_0x^2_1+x^2_0x_2^2+x^2_0x^3_3+x^2_1x^2\\
~&~&+x^2_1x_3^2+x_2^2x_3^2)+8x_0x_1x_2x_3.
\end{array}
\end{align}
It follows that
\begin{align}\label{4.3.4}
\begin{array}{l}
{\cal Z}_{{\cal A}_8}=\left\{(x_0,x_1,x_2,x_3)^T\in \R^4\;|\;\right.\\
\left. \xi_{{\cal A}_8}(x_0,x_1,x_2,x_3)=0\right\}.
\end{array}
\end{align}

\end{itemize}

\end{exa}

\subsection{Some Other Examples}

To see there are also  ${\cal A}\in \H_n$, for $n>4$, such examples are presented as follows.
First example is a set of simplest PHAs, which are called trivial PHAs.

\begin{exa}\label{e4.4.1} Define an $n+1$ dimensional algebra ${\cal A}_{n+1}^0$ as follows:
Let $\i_k$, $k=1,2,\cdots,n$ be its hyperimaginary units. Set
$$
\i_s*\i_t=0,\quad s,t=1,2,\cdots,n.
$$
Then it is easy to verify that ${\cal A}_{n+1}^0\in \H_{n+1}$. Moreover, its PSM, $P_{{|cal A}_{n+1}^0}$ can be determined by the following:
\begin{align}\label{4.4.1}
\Col_i\left(P_{{\cal A}_{n+1}^0}\right)=
\begin{cases}
\d_{n+1}^i,\quad i=1,2,\cdots,n+1;\\
\d_{n+1}^{r+1},\quad i=r(n+1)+1,\;r=1,2,\cdots,n,\\
0,\quad \mbox{Otherwise}.
\end{cases}
\end{align}
Its characteristic function is
\begin{align}\label{4.4.2}
\xi_{{\cal A}_D^n}=x_0^{n+1}.
\end{align}
Hence,
\begin{align}\label{4.4.3}
{\cal Z}_{{\cal A}_D^n}=\{x_0+x_1\i_1+\cdots+x_n\i_n\;|\;x_0=0\}.
\end{align}
If $x\in {\cal Z}^c_{{\cal A}_D^n}$, say
$$
x=x_0+x_1\i_1+\cdots+x_n\i_n,\quad x_0\neq 0,
$$
then
$$
x^{-1}=\frac{1}{x_0}-\dsum_{i=1}^n\frac{x_i}{x_0^2}\i_i.
$$
\end{exa}

Next, we give an example for $n=5$.

\begin{exa}\label{e4.4.2} Consider a hypercomplex algebra ${\cal A}$, which has its PSM as
\begin{align}\label{4.4.4}
\begin{array}{ccl}
P_{{\cal A}}&:=&\d_5[1,2,3,4,5,2,0,0,3,0,3,0,\\
&&~~~0,0,0,4,3,0,0,0,5,0,0,0,0],\\
\end{array}
\end{align}
where $\d_5^0={\bf 0}_5$.

A straightforward computation shows that ${\cal A}$ is commutative and associative, hence ${\cal A}\in \H_5$.
Then it is easy to calculate that
\begin{align}\label{4.4.5}
\xi_{{\cal A}}(x)=x_0^3(x_0^2-x_1x_3).
\end{align}
So its zero set is
\begin{align}\label{4.4.6}
{\cal Z}_{{\cal A}}=\{(x_0,x_1,x_2,x_3,x_4)^T\in \R^5\;|\;x_0=0, ~\mbox{or}~x_0^2=x_1x_3\}.
\end{align}
\end{exa}

\section{Conclusion}

In this paper the perfect hypercomplex algebra is considered. Using STP, necessary and sufficient conditions on the structure matrix of an algebra to be a PHA is proposed. Based on the matrix expression of homomorphisms between algebras, certain lower dimensional PHAs are classified up to isomorphism. Their characteristic functions and zero sets are discussed.

\end{document}